\documentclass[12pt]{article}
\usepackage{graphicx}
\usepackage{psfrag}

\usepackage{amsmath,amssymb}
\usepackage{mathrsfs}
\usepackage{amsfonts}
\usepackage{latexsym}
\usepackage{amsthm}
\usepackage{pictex}

\usepackage[T2A]{fontenc}     
\usepackage[cp1251]{inputenc} 
\usepackage{graphicx}

\newcommand{\er}[1]{{\rm(\ref{#1})}}
\def\lb{\label}
\theoremstyle{plain}
\newtheorem{theorem}{\bf Theorem}[section]
\newtheorem{lemma}[theorem]{\bf Lemma}
\theoremstyle{remark}

\begin{document}

\def\a{\alpha} \def\cA{{\cal A}} \def\bA{{\bf A}}  \def\mA{{\mathscr A}}
\def\b{\beta}  \def\cB{{\cal B}} \def\bB{{\bf B}}  \def\mB{{\mathscr B}}
\def\g{\gamma} \def\cC{{\cal C}} \def\bC{{\bf C}}  \def\mC{{\mathscr C}}
\def\G{\Gamma} \def\cD{{\cal D}} \def\bD{{\bf D}}  \def\mD{{\mathscr D}}
\def\d{\delta} \def\cE{{\cal E}} \def\bE{{\bf E}}  \def\mE{{\mathscr E}}
\def\D{\Delta} \def\cF{{\cal F}} \def\bF{{\bf F}}  \def\mF{{\mathscr F}}
\def\c{\chi}   \def\cG{{\cal G}} \def\bG{{\bf G}}  \def\mG{{\mathscr G}}
\def\z{\zeta}  \def\cH{{\cal H}} \def\bH{{\bf H}}  \def\mH{{\mathscr H}}
\def\e{\eta}   \def\cI{{\cal I}} \def\bI{{\bf I}}  \def\mI{{\mathscr I}}
\def\p{\psi}   \def\cJ{{\cal J}} \def\bJ{{\bf J}}  \def\mJ{{\mathscr J}}
\def\vT{\Theta}\def\cK{{\cal K}} \def\bK{{\bf K}}  \def\mK{{\mathscr K}}
\def\k{\kappa} \def\cL{{\cal L}} \def\bL{{\bf L}}  \def\mL{{\mathscr L}}
\def\l{\lambda}\def\cM{{\cal M}} \def\bM{{\bf M}}  \def\mM{{\mathscr M}}
\def\L{\Lambda}\def\cN{{\cal N}} \def\bN{{\bf N}}  \def\mN{{\mathscr N}}
\def\m{\mu}    \def\cO{{\cal O}} \def\bO{{\bf O}}  \def\mO{{\mathscr O}}
\def\n{\nu}    \def\cP{{\cal P}} \def\bP{{\bf P}}  \def\mP{{\mathscr P}}
\def\r{\rho}   \def\cQ{{\cal Q}} \def\bQ{{\bf Q}}  \def\mQ{{\mathscr Q}}
\def\s{\sigma} \def\cR{{\cal R}} \def\bR{{\bf R}}  \def\mR{{\mathscr R}}
\def\S{\Sigma} \def\cS{{\cal S}} \def\bS{{\bf S}}  \def\mS{{\mathscr S}}
\def\t{\tau}   \def\cT{{\cal T}} \def\bT{{\bf T}}  \def\mT{{\mathscr T}}
\def\f{\phi}   \def\cU{{\cal U}} \def\bU{{\bf U}}  \def\mU{{\mathscr U}}
\def\F{\Phi}   \def\cV{{\cal V}} \def\bV{{\bf V}}  \def\mV{{\mathscr V}}
\def\P{\Psi}   \def\cW{{\cal W}} \def\bW{{\bf W}}  \def\mW{{\mathscr W}}
\def\o{\omega} \def\cX{{\cal X}} \def\bX{{\bf X}}  \def\mX{{\mathscr X}}
\def\x{\xi}    \def\cY{{\cal Y}} \def\bY{{\bf Y}}  \def\mY{{\mathscr Y}}
\def\X{\Xi}    \def\cZ{{\cal Z}} \def\bZ{{\bf Z}}  \def\mZ{{\mathscr Z}}
\def\O{\Omega}
\def\ve{\varepsilon}
\def\vt{\vartheta}
\def\vp{\varphi}
\def\vk{\varkappa}

\def\mM{M}
\def\mB{B}
\def\mR{R}

\def\mA{{\mathscr A}}
\def\mB{{\mathscr B}}
\def\mC{{\mathscr C}}
\def\mD{{\mathscr D}}
\def\mE{{\mathscr E}}
\def\mF{{\mathscr F}}
\def\mG{{\mathscr G}}
\def\mH{{\mathscr H}}
\def\mI{{\mathscr I}}
\def\mJ{{\mathscr J}}
\def\mK{{\mathscr K}}
\def\mL{{\mathscr L}}
\def\mM{{\mathscr M}}
\def\mN{{\mathscr N}}
\def\mO{{\mathscr O}}
\def\mP{{\mathscr P}}
\def\mQ{{\mathscr Q}}
\def\mR{{\mathscr R}}
\def\mS{{\mathscr S}}
\def\mT{{\mathscr T}}
\def\mU{{\mathscr U}}
\def\mV{{\mathscr V}}
\def\mW{{\mathscr W}}
\def\mX{{\mathscr X}}
\def\mY{{\mathscr Y}}
\def\mZ{{\mathscr Z}}

\def\Z{{\Bbb Z}}
\def\R{{\Bbb R}}
\def\C{{\Bbb C}}
\def\T{{\Bbb T}}
\def\N{{\Bbb N}}
\def\S{{\Bbb S}}
\def\H{{\Bbb H}}
\def\J{{\Bbb J}}

\def\qqq{\qquad}
\def\qq{\quad}
\newcommand{\ma}{\begin{pmatrix}}
\newcommand{\am}{\end{pmatrix}}
\newcommand{\ca}{\begin{cases}}
\newcommand{\ac}{\end{cases}}
\let\ge\geqslant
\let\le\leqslant
\let\geq\geqslant
\let\leq\leqslant
\def\ma{\left(\begin{array}{cc}}
\def\am{\end{array}\right)}
\def\iint{\int\!\!\!\int}
\def\lt{\biggl}
\def\rt{\biggr}
\let\geq\geqslant
\let\leq\leqslant
\def\[{\begin{equation}}
\def\]{\end{equation}}
\def\wh{\widehat}
\def\wt{\widetilde}
\def\pa{\partial}
\def\sm{\setminus}
\def\es{\emptyset}
\def\no{\noindent}
\def\ol{\overline}
\def\iy{\infty}
\def\ev{\equiv}
\def\/{\over}
\def\ts{\times}
\def\os{\oplus}
\def\ss{\subset}
\def\h{\hat}
\def\Re{\mathop{\rm Re}\nolimits}
\def\Im{\mathop{\rm Im}\nolimits}
\def\supp{\mathop{\rm supp}\nolimits}
\def\sign{\mathop{\rm sign}\nolimits}
\def\Ran{\mathop{\rm Ran}\nolimits}
\def\Ker{\mathop{\rm Ker}\nolimits}
\def\Tr{\mathop{\rm Tr}\nolimits}
\def\const{\mathop{\rm const}\nolimits}
\def\dist{\mathop{\rm dist}\nolimits}
\def\diag{\mathop{\rm diag}\nolimits}
\def\Wr{\mathop{\rm Wr}\nolimits}
\def\BBox{\hspace{1mm}\vrule height6pt width5.5pt depth0pt \hspace{6pt}}

\def\Diag{\mathop{\rm Diag}\nolimits}

\def\Twelve{
\font\Tenmsa=msam10 scaled 1200 \font\Sevenmsa=msam7 scaled 1200
\font\Fivemsa=msam5 scaled 1200 \textfont\msbfam=\Tenmsb
\scriptfont\msbfam=\Sevenmsb \scriptscriptfont\msbfam=\Fivemsb

\font\Teneufm=eufm10 scaled 1200 \font\Seveneufm=eufm7 scaled 1200
\font\Fiveeufm=eufm5 scaled 1200
\textfont\eufmfam=\Teneufm \scriptfont\eufmfam=\Seveneufm
\scriptscriptfont\eufmfam=\Fiveeufm}

\def\Ten{
\textfont\msafam=\tenmsa \scriptfont\msafam=\sevenmsa
\scriptscriptfont\msafam=\fivemsa

\textfont\msbfam=\tenmsb \scriptfont\msbfam=\sevenmsb
\scriptscriptfont\msbfam=\fivemsb

\textfont\eufmfam=\teneufm \scriptfont\eufmfam=\seveneufm
\scriptscriptfont\eufmfam=\fiveeufm}

\title {Nanoribbons in external electric fields}

\author{
 Evgeny Korotyaev
\begin{footnote}
{School of Mathematics, Cardiff Univ., Senghennydd Road, Cardiff,
CF24 4AG, UK, e-mail: korotyaeve@cf.ac.uk}
\end{footnote}
 \and Anton Kutsenko
\begin{footnote}
{ Department of
 Mathematics of Sankt-Petersburg State University, Russia e-mail: kucenkoa@rambler.ru}
\end{footnote}
}

\maketitle

\begin{abstract}
\no    We consider the Schr\"odinger operator on nanoribbons
 (tight-binding models) in an external electric potentials  $V$.
 The corresponding electric field is perpendicular to the axis of the  nanoribbon.
If $V=0$, then the spectrum of the Schr\"odinger operator consists
of two spectral bands and the flat band (i.e., the eigenvalue with
infinite multiplicity) between them. If we switch on an weak
electric potential $V\to 0$, then we determine the asymptotics of
the spectral bands for small fields. In particular, we describe all
potentials when the unperturbed eigenvalue remains the flat band and
when one becomes the small band of the continuous spectrum.
\end{abstract}


\section {Introduction}
\setcounter{equation}{0}

After their discovery \cite{Ii}, carbon nanotubes remain in both
theoretical and applied research \cite{SDD}. Last years physicists
consider also nanoribbons and in particular nanoribbons in external
electric fields, see \cite{SCL}. We consider the Schr\"odinger
operator $H=\D+V$ on the nanoribbon $\G$ (a tight-binding model of
single-wall nanoribbons, see \cite{SDD}, \cite{N}) in an external
electric potential $V$. The electric field is perpendicular to the
axis of the nanoribbon.  Our model nanoribbon $\G\ss\R^2$ is a
graph, which is a set of vertices $\vk_{n,k}$ and  bonds (edges)
$\G_{n,k,j}$  given by
\[
 \ca\vk_{n,2k+1}=(\sqrt3(2n+k),3k), & k\in\N^0_N\\
 \vk_{n,2k}=(\sqrt3(2n+k),3k-2), & k\in\N_{N}\ac,\ \ \ca
 \G_{n,k,1}=[\vk_{n,2k},\vk_{n,2k+1}]\\\G_{n,k,2}=[\vk_{n,2k},\vk_{n,2k-1}]\\
 \G_{n,k,3}=[\vk_{n,2k},\vk_{n+1,2k-1}]
 \ac,\ \ n\in\Z,
\]
where  $\N_k=\{1,..,k\}\ss\N$ and $\N^0_k=\N_k\cup\{0\}$. See Fig
\ref{f001} for the case $N=3$. We define the discrete Hilbert space
$\ell^2(\G)$ consisting of functions $f_{n,k}$ on the set of
vertices  $\vk=\{\vk_{n,k}, (n,k)\in \Z\ts \N_p\}, p=2N+1$ equipped
with the norm $\|f\|_{\ell^2(\G)}^2=\sum |f_{n,k}|^2$.
\begin{figure}[t]
\lb{f001}
         \centering\includegraphics[clip]{ribbon.0}
         \begin{center}
         {\small A  horizontal nanoribbon at $N=3, k\in \{1,2,..,7\}$ in the
         vertical electric field. }
         \end{center}
         \end{figure}
The Laplace operator $\D$ and the potential $V$ are given by
\begin{multline}
 \lb{H}
 (\D f)_{n,2k+1}=f_{n,2k}+f_{n-1,2k+2}+f_{n,2k+2},
 \qq k\in \N^0_{N},\\
 (\D f)_{n,2k}=f_{n,2k-1}+f_{n+1,2k-1}+f_{n,2k+1},
 \ \ k\in\N_N,
\end{multline}
\[
 \lb{dbc}
 f_{n,0}=f_{n,p+1}=0,\qq
 p=2N+1, \qq  \ n\in \Z, \qq f=(f_{n,k})_{(n,k)\in \Z\ts
\N_p}\in\ell^2(\G),
\]
\[
\lb{V} (Vf)_{n,k}=v_k f_{n,k},\qqq (n,k)\in \Z\ts\N_{p},
\]
Our electric potential is given by \er{V}, since the electric field
is perpendicular to the axis of the nanoribbon. In fact we consider
the Schr\"odinger operator $H$  on the set $\Z\ts\N_{p}$ with the
Dirichlet boundary conditions \er{dbc}. We formulate our preliminary
result.

\begin{theorem}
\label{T1} i)  The operator $H=\D+V$ is unitarily equivalent to the
operator $\int_{[0,2\pi)}^{\os}J_a{dt\/2\pi }$, where $p\ts p$
matrix $J_a, a=2|\cos {t\/2}|$  is a Jacobi operator, acting on
$\C^{p}$ and given by
\begin{multline}
 \lb{Jk}
 (J_a y)_n=a_{n-1}y_{n-1}+a_{n}y_{n+1}+v_ny_n,\qq
 y=(y_n)_1^{p}\in\C^{p},\qq y_0=0=y_{p+1},\qq p=2N+1,\\
a_1=a_{2n+1}=a=2|\cos {t\/2}|,
  \qq  \qq a_{2n}=1,\ \  \  n\in\N_{N}.
\end{multline}
\no  ii) The spectrum of $H$ is given by
\[
\lb{spp}
 \s(H)=\bigcup_{k=-N}^{N}\s_n,\qqq     \s_k=\l_k([0,2])=\ca [ \l_k^{-},\l_k^{+}],
  & k\ge0\\
 [\l_k^{+},\l_k^{-}], & k<0\ac,
\]
where $\l_{-N}(a)\le \l_{-N+1}(a)\le ...\le \l_{N}(a)$ are
eigenvalues of $J_a$.  Moreover,  $\l_{n}(\cdot)$ is real analytic
in $a\in [0,2]$ and if $a\ne 0$, then $\l_{-N}(a)< \l_{-N+1}(a)<
...< \l_{N}(a)$.
\end{theorem}

{\no\bf Remark.} 1)  The matrix of the operator $J_a$ is given by
\[
\lb{Jk1} J_a=\left(\begin{array}{ccccc}
                                    v_1 & a & 0 & .. & 0 \\
                                    a & v_2 & 1 & .. & 0\\
                                    0 & 1 & v_3 & .. & 0\\
                                    .. & .. & .. & .. & ..\\
                                    0 & .. & 0 & 1 & v_p
                 \end{array}\right)=J_a^0+\diag(v_n)_1^p,
                 \qq v=(v_n)_1^p\in \R^p.
\]
2) Exner \cite{Ex} obtained a duality between Schr\"odinger
 operators on graphs and certain (depending on energy) Jacobi matrices. In our case the
  Jacobi matrices do not depend on energy.

We describe the spectrum of the unperturbed operator $\D$. Let
$\l_{-N}^0(a)\le \l_{-N+1}^0(a)\le ...\le \l_{N}^0(a)$ be
eigenvalues of $J_a^0$ corresponding to $\D$  and let
$\Z_{N}=\{-N,...,N\}$.

\begin{theorem}
\label{T2}
 Let $V=0$ and let $c_k=\cos\frac{k\pi}{N+1}$,
$s_k=\sin\frac{k\pi}{N+1}, k\in \Z_N$. Then
\[
\lb{T2-1}
 \l_{-k}^0(a)=\l_k^0(a)=(a^2-2ac_k+1)^{1\/2},\ \
 k\in \N_N,\ \ \ \ \l_{0}^0(a)=0, \qqq all \qq a\in [0,2],
\]
\begin{multline}
\lb{T2-2}
 \s(\D)=\s_{ac}(\D)\cup \s_{pp}(\D)=
 \bigcup_{-N}^{N} \s_k^0, \qqq \s_{pp}(\D)=\s_0^0=\{0\},
 \\
 \s_{ac}(\D)=[-\l_N^0(2),\l_N^0(2)]\sm (-s_1,s_1),\qqq
 \l_N^0(2)=(5+4c_1)^{1\/2},
\end{multline}
\[
\lb{T2-3} \s_{-k}^0=-\s_k^0,   \qq
\s_k^0=[\l_k^{0,-},\l_k^{0,+}]=\ca [s_k, \l_k^0(2)],  &  if \ \ c_k<0\\
            [1,\l_k^0(2)], &  if \ c_k\ge0 \ac,\ \qqq all \qq k\in\N_N.
\]
\end{theorem}

\no{\bf Remark.} 1) The a.c. spectrum $\s_{ac}(\D)$ consists of two
spectral bands  $[-\l_N^0(2), -s_1]$ and $[s_1,\l_N^0(2)]$ separated
by the gap $\g^0=(-s_1,s_1)$ and $\D$ has the flat band
$\s_0^0=\{0\}$.

\no 2) Note that the gap length $|\g^0|=2s_1={2\pi \/N}(1+o(1))$ as
$N\to \iy$ (see also \cite{SCL}).

\no 3) The following spectral intervals $\s_N^0\sm
\s_{N-1}^0=(\l_{N-1}^0(2),\l_N^0(2)]$ and $\s_1^0\sm
\s_{2}^0=[s_1,s_2)$ has the spectrum of $\D$ of multiplicity 2.

Consider the spectrum of $H$ for a small potential $V\to 0$. In this
case the spectrum of  $H$ will be the small perturbation  of the
spectrum $\D$. Under the perturbation $V$ all spectral  bands
$\s_k^0, k\ne 0$ of $\D$ will be again spectral bands $\s_k, k\ne 0$
of the operator $H$. The more complicated  case of $\s_0$ will be
considered in Theorem \ref{T4a} and \ref{T4}.  In the following
theorem  we determine asymptotics of $\s_k$. Below we will {\bf
sometimes write} $\l_k(a,v), \s_k(v),..$, instead of $\l_k(a), \s_k,
..$, when several potentials are being dealt with. Let
$\|v\|^2=\sum_1^pv_n^2$.

\begin{theorem}
\label{T3}
 i)  Each eigenvalue $\l_k(a,v), \pm k \in \N_n$ of
$J_a$  is a real analytic function of $(a,v)\in (0,2]\ts \R^p$ and
$\pm \pa_a^2\l_k(a,v)>0, \pm k \in \N_n$ for $a\in[\d,2]$ for each
small $\d>0$ and sufficiently small $v$. Moreover,
$\l_0(\cdot,\cdot)$ is analytic in the neighborhood of $(0,0)$.

\no ii) If $\|v\|^2=\sum_1^p|v_n|^2 \to 0$ and $k\ne0$, then
\[
 \lb{T3-1}
 \l_k^-(v)=s_k\sign k +{\sum_{n=1}^{N+1}(c^2_{nk}v_{2n-1}+
 s^2_{nk}v_{2n})\/N+1}+O(\|v\|^2),\qq
 |k|<\frac{N+1}2,
\]
\[
 \lb{T3-2}
 \l_k^-(v)=\sign k+O(\|v\|),\ |k|\ge\frac{N+1}2,
\]
\[
 \lb{T3-3}
 \l_k^+(v)=(5-4c_k)^{\frac12}\sign
 k+\sum_{n=1}^{p}{\c_nv_n\/N+1}+O(\|v\|^2),\ \
 \ca
 \c_{2n}=s_{kn}^2 \\
 \c_{2n+1}={(s_{nk}-2s_{(n+1)k})^2\/(5-4c_k)}.
 \ac
\]
\end{theorem}

 The following Theorems \ref{T4a} and \ref{T4} are our main result
about the spectral interval $\s_0\ne \es$ (or the flat band $\s_0$)
as $v\to 0$. Let $\cH=\{f\in\ell^2(\G): Hf=v_1f\}$.

\begin{theorem}
\lb{T4a}
i) The spectral band $\s_{0}(v)$ of $H$ is a flat band and
$\s_{0}(v)=\{\l_{0}(a)\}=\{v_1\}$ (an eigenvalue of $H$) iff
$v_{2k+1}=v_1$ for all $k\in\N_N$.

\no ii) Let $\s_{0}(v)=\{v_1\}$ and let $Sh=(h_{n+1})_{n\in\Z}$ for
$h=(h_{n+1})_{n\in\Z}\in \ell^2(\Z)$. Then
\begin{multline}
\cH=\{f\in\ell^2(\G): Hf=v_1f\}= \{f=(f_{n,k})_{(n,k)\in \Z\ts
\N_p}\in\ell^2(\G): \\ (f_{n,2k+1})_{n\in\Z}=(-I-S)^kh,\qq
(f_{n,2k})_{n\in\Z}=0,\ k\in\N_N, \ h=(h_{n})_{n\in\Z}\in\ell^2(\Z)
\}.
\end{multline}

\no iii)  Let $\s_{0}(v)=\{v_1\}$. Define compactly supported
functions $\p^m=(\p_{n,k}^m)_{(n,k)\in \Z\ts \N_p}\in\cH, m\in \Z$
 by
\[
\lb{T4a-1} \p_k^m=(\p_{n,k}^m)_{n\in \Z},\qqq \p_{2k}^m=0, \qq
\p_{2k+1}^m=(-I-S)^ke_m,\qq k\in \N_N^0,
\]
where $e_m=(\d_{n,m})_{n\in\Z}\in\ell^2(\Z)$. Then each $f\in \cH$
has the form
\[
\lb{T4a-2}
f=\sum_{m\in\Z}\wh f_m\p^{m},\qqq \wh f_{m}= f_{m,1}.
\]
 Moreover, the mapping $f\to (\wh
f_{m})_{m\in\Z}$ is a linear isomorphism between $\cH$ and
$\ell^2(\Z)$.

\end{theorem}

\no {\bf Remark.} Due to \er{T4a-2} each compactly supported
eigenfunction has the form  $f=\sum_{\a}^\b\wh f_m\p^{m}$ for some
$\a, \b$.

\begin{theorem}
\label{T4} i) Let $\|v\|^2=\sum_1^p|v_n|^2\to 0$. Then
\[
\lb{T5-1}
\l_{0}(a,v)=F(a,v)+O(\|v\|^2),\qqq
 F(a,v)={\sum_{k=0}^Nv_{2k+1}a^{2k}\/\sum_{k=0}^Na^{2k}},
\]
\[ \lb{T5-1a}
\l_0^-(v)=\min_{a\in[0,2]}F(a,v)+O(\|v\|^2),\qq
\l_0^+(v)=\max_{a\in[0,2]}F(a,v)+O(\|v\|^2).
\]
\no ii) Let $0=v_{1}\le v_{3}\le..\le v_{p}$ and  $0< v_{p}$. Then
the following asymptotics hold true
\[
\lb{T5-2}
\l_0^-(v)=O(v_p^2),\qqq
 \l_{0}^+(v)={3\/4^{N+1}-1}\sum_{k=1}^N4^{k}v_{2k+1}+O(v_p^2) \qq
 as \qqq v\to 0.
\]
Moreover, if $v_{1}<v_{3}\le v_p\le Cv_3$ for some $C>0$, then the
function $\l_{0}(\cdot,v)$ is strongly increasing on $[0,2]$ for
sufficiently small $v_p$ and the spectral band $\s_0(v)$ has the
spectrum of $H$ of multiplicity 2.
\end{theorem}

{\bf Remark.} 1) The first term of  asymptotics of $\s_0$ depends on
the odd components $v_{2k+1}$ of the potential. Moreover, $\s_0(v)$
is a flat band iff all odd components $v_{2k+1}=0$.

\no 2) If $v_1\not=0$ and all odd $v_{2k+1}=0, k\ge 1$, then
$\s^0(v)=|v_1|(1-\frac{3}{4^{N+1}-1})+O(\|v\|^2)$.

{\bf Example of the constant electric fields.} Consider the electric
field given by  $\cE=\cE_0(0,1)\in \R^2$, where $\cE_0>0$ is the
amplitude. The corresponding electric potential $\cV(x,y)=\cE_0 y,
(x,y)\in \R^2$ is increasing in the vertical direction and
$$
v_{2k+1}=\ve k,\qq k=0,1,.,N,
$$
for some $\ve>0$. We assume that the electric field is weak  $V\to
0$  and then $\ve \to0$. Thus asymptotics \er{T5-2} give
\[
\lb{cef} \l_0^-(v)=O(\ve^2),\qqq
 \l_{0}^+(v)=4\ve C_p+O(\ve^2),\qq C_p={4N(4^{N}-1)-3\/3(2^{p+1}-1)}
\]
as $\ve\to 0$. Thus if we switch on the constant electric field,
then $\s(H)=\s_{ac}(H)$ and $|\s_0(v)|=4\ve C_p+O(\ve^2)>0$ for
small $\ve >0$.

{\bf Strong electric fields.} We consider the nanoribbon in strong
electric fields. Our operator has the form $H(t)=\D+tV$ as the
coupling constant $t\to\iy$. In this case the corresponding Jacobi
operator depends on $t$ and is given by
\[
\lb{jt} (J_a(t)y)_n=a_{n-1}y_{n-1}+a_{n}y_{n+1}+tv_ny_n,\
y=(y_n)_{n\in \Z}\in \C^p.
\]
We study how the spectral bands
$\s_{k}(tv)=[\l_{k}^+(t),\l_{k}^{-}(t)]$, $k\in\Z_N$ of the operator
$H(t)$ depend on the coupling constant $t\to \iy$.

\begin{theorem}
\lb{T5} Let $H(t)=\D+tV$ with $v_1<..<v_p$ and let $t\to\iy$. Then
\begin{multline}
\lb{T5.1}
 \l_{k-N-1}^{\pm s}(t)=t v_k-{\x_k^{\pm}+O(t^{-1})\/t},\qqq
 \x_k^{\pm}={r_{k}^{\pm}\/v_{k-1}- v_k}+{r_{k+1}^{\pm}\/v_{k+1} -
 v_k},\ \  s=(-1)^k,\\
v_0=v_{p+1}=0,\qq  r_1^{\pm}=r_{p+1}^{\pm}=0,\ \ r_{2n+1}^{\pm}=1,\
\
 r_{2n}^-=0,\ \ r_{2n}^+=4,\ \ n\in\N_{N},
\end{multline}
\[ \lb{T5.3}
 |\s_{k-N-1}(tv)|={4+O(t^{-1})\/t|v_{k-(-1)^k}-v_k|},\ \ k\in
 \N_p\sm\{p\},\qq and \qq
 |\s_{N}(tv)|=O(t^{-2})
\]
for any $k\in\N_p$, where each spectral band $\s_k(t)$, $k\ne N$ has
multiplicity $2$.
\end{theorem}

Note that \er{T5.1} implies $\s_k(t)\cap\s_j(t)=\es$ for any $k\ne
j$, $k,j\in\Z_N$ as $t\to \iy$.

A carbon nanoribbon is a strip of the honeycomb lattice, i.e., a cut
graphene sheet, see a horizontal ribbon in Fig. 1. There are
mathematical results about Schr\"odinger operators on  carbon
nanotubes (zigzag, armchair and chiral) of Korotyaev and Lobanov
 \cite{KL}, \cite{KL1}, \cite{K1}, Kuchment and Post \cite{KuP} and
Pankrashkin \cite{Pk}. All these papers consider so-called
continuous models. But in the physical literature the most commonly
used model is the tight-binding model. For applications of our model
see ref. in \cite{ARZ}, \cite{H}, \cite{SDD}.

In the proof of our theorems we determine various asymptotics for
periodic Jacobi operators with specific coefficients, see \er{Jk}.
Note that there exist a lot of papers devoted to asymptotics and
estimates both for periodic Jacobi operators and Schr\"odinger
operators see e.g. \cite{KKu1}, \cite{KKr}, \cite{La}, \cite{vMou}.

We present the plan of our paper. In Sect. 2 we prove Theorems
\ref{T1} and \ref{T2}. In the proof Theorem \ref{T1} we use
arguments from \cite{KL}. In Sect. 3 we prove Theorems \ref{T3}
-\ref{T5}.

\section {Preliminaries}
\setcounter{equation}{0}

{\bf Proof of Theorem \ref{T1}.} i) 
For each $(f_{n,k})_{(n,k)\in \Z\ts \N_p}\in \ell^2(\G)$ we
introduce the function $\p_k=(f_{n,k})_{n\in\Z}\in\ell^2(\Z)$,
$k\in\N_{p}, p=2N+1 $ and $\p=(\p_k)_{k\in\N_p}\in(\ell^2(\Z))^p$.
Using \er{H} and $(Vf)_{n,k}=v_kf_{n,k}$ for any
$(n,k)\in\Z\ts\N_p$, we obtain that the operator
$H:\ell^2(\G)\to\ell^2(\G)$ is unitarily equivalent to the operator
$K:\ell^2(\Z)^p\to\ell^2(\Z)^p$, given by
$$
 (K\p)_{2k+1}=\p_{2k}+v_{2k+1}\p_{2k+1}+(I+S^{-1})\p_{2k+2},\qq \p_0=\p_{p+1}=0,
 \qq k\in \N^0_{N},
$$
$$
 (K\p)_{2k}=(I+S)\p_{2k-1}+v_{2k}\p_{2k}+\p_{2k+1},
 \ \ k\in\N_N.
$$
where the operators $Ih=h$ and $Sh=(h_{n+1})_{n\in\Z}$,
$h=(h_n)_{n\in\Z}$. We rewrite $K$ in the matrix form by
\[
 \lb{Km}
 K(\p_k)_{1}^{p}=\left(\begin{array}{ccccc}
                                    v_1 & I+S^{-1} & 0 & .. & 0 \\
                                    I+S & v_2 & I & .. & 0\\
                                    0 & I & v_3 & .. & 0\\
                                    .. & .. & .. & .. & ..\\
                                    0 & .. & 0 & I & v_{p}
                 \end{array}\right)
                 \left(\begin{array}{c}\p_1\\ \p_2\\ \p_3\\..\\ \p_{p}\end{array}\right),\ \
                 \p_k\in\ell^2(\Z).
\]
Note  that $K^*=K$, since $S^*=S^{-1}$. Introduce the unitary
operator  $\F:\ell^2(\Z)^p\to \int_{[0,2\pi)}^{\os}\mH_0
{dt\/2\pi},\ \ \mH_0=\C^p$, by $\F (\p_k)_{1}^p=(\f \p_k)_{1}^p$,
where  $\f:\ell^2(\Z)\to L^2(0,2\pi)$ is an unitary operator given
by
$$
\f  h=\sum_{n\in\Z}h_n{e^{int}\/\sqrt{2\pi}}, \qqq h=
(h_n)_{n\in\Z}\in\ell^2(\Z),\  t\in[0,2\pi].
$$
Then we deduce that
$
\F_pK\F_p^{-1}=\int_{[0,2\pi)}^{\os}\wt J_t{dt\/2\pi},
$
where the operator $\wt J_t: \mH_0\to \mH_0$ has the matrix given by
$$
\wt J_t y= \left(\begin{array}{ccccc}
                                    v_1 & 1+e^{it} & 0 & .. & 0 \\
                                    1+e^{-it} & v_2 & 1 & .. & 0\\
                                    0 & 1 & v_3 & .. & 0\\
                                    .. & .. & .. & .. & ..\\
                                    0 & .. & 0 & 1 & v_{p}
                 \end{array}\right)
                 \left(\begin{array}{c}y_1\\ y_2\\ y_3\\..\\ y_{p}\end{array}\right),
                 \qq y=(y_k)_1^{p}\in \C^{p}.
$$
The matrix $\wt J_t$ is unitarily equivalent to the matrix $J_a$,
$a=2|\cos\frac t2|=|1+e^{it}|$, given by
\[
\lb{ja0}
 J_a=\left(\begin{array}{ccccc}
                                    v_1 & a & 0 & .. & 0 \\
                                    a & v_2 & 1 & .. & 0\\
                                    0 & 1 & v_3 & .. & 0\\
                                    .. & .. & .. & .. & ..\\
                                    0 & .. & 0 & 1 & v_{p}
                 \end{array}\right)=J_a^0+\diag(v_n)_{1}^{p},\qq
                 J_a^0=\left(\begin{array}{ccccc}
                                    0 & a & 0 & .. & 0 \\
                                    a & 0 & 1 & .. & 0\\
                                    0 & 1 & 0 & .. & 0\\
                                    .. & .. & .. & .. & ..\\
                                    0 & .. & 0 & 1 & 0
                 \end{array}\right).
\]
Thus we deduce that the operator $H$ is unitarily equivalent to the
operator $\int_{[0,2\pi)}^{\os}J_a{dt\/2\pi}$.

ii) From the spectral theory of Jacobi operators \cite{vM} we obtain
that eigenvalues $\l_k(a), k\in \Z_N$ of $J_a$ satisfy
$\l_{-N}(a)\le \l_{-N+1}(a)\le ...\le \l_{N}(a)$. Moreover,  if
$a\ne 0$, then $\l_{-N}(a)< \l_{-N+1}(a)< ...< \l_{N}(a)$. The
perturbation theory gives that $\l_{n}(\cdot)$ is real analytic in
$a\in [0,2]$. Then the standard arguments for periodic operator (see
Sect. XIII, 16 \cite{RS}) and i) give the proof of ii). \BBox

Below we need the well known facts about Jacobi matrices. Let
$a\not=0$. Introduce the fundamental solutions $\vp_n=\vp_n(\l,a,v),
\vt_n=\vt_n(\l,a,v)$, $(\l,a,v)\in\C\ts \R_+\ts \R^{2N+1}$ of the
system of equations
$$
 a y_{2k}-(\l-v_{2k-1})y_{2k-1}+y_{2k-2}=0,\qq
 y_{2k+1}-(\l-v_{2k})y_{2k}+ay_{2k-1}=0,\ \ k\in\N_{N+1},
$$
with conditions: $ \vt_0=0,\ \ \vt_1=1,  \vp_0=0,\ \ \vp_1=1,\ \ $
and $v_{p+1}=v_1$. The solution $\vp_k$ has the following  form
\[
\lb{mo}
 \left(\begin{array}{c} \vp_{2k} \\
 \vp_{2k+1}\end{array}\right)=M_k\left(\begin{array}{c} 0\\
 1\end{array}\right),
 \]
where  $M_k$ is the monodromy matrix, given by
\[
\lb{mm1}
 M_k=\ma\vt_{2k}&\vp_{2k}\\
 \vt_{2k+1}&\vp_{2k+1}\am=T_k..T_1,\ \
  T_k={1\/a}\ma -1& \l-v_{2k-1}\\
 v_{2k}-\l& (\l-v_{2k})(\l-v_{2k-1})-a^2
 \am.
\]

Consider the spectrum for the unperturbed operator $\D$.

\no {\bf Proof of Theorem \ref{T2}.} Using \er{mm1}, \er{mo}, we
have that
$$
 \left(\begin{array}{c} \vp_{2k} \\
 \vp_{2k+1}\end{array}\right)=T_1^{k}\left(\begin{array}{c} 0\\
 1\end{array}\right),\ \
 T_1={1\/a}\ma -1& \l\\  -\l& \l^2-a^2
 \am.
$$
The eigenvalues $\t_{\pm}$ and eigenvectors $e_{\pm}$ of the matrix
$M_1$ are given by
\[
 \lb{t0001}
 \t_{\pm}=\frac{r\pm\sqrt{r^2-4}}{2},\ \
 e_{\pm}=\left(\begin{array}{c} \frac{\l}a \\
 \frac1a+\t_{\pm}\end{array}\right),\ \ r=\Tr M_1=\frac{\l^2}a-a-\frac1a
\]
and satisfy
$$
\t_+\t_-=\det M_1=1,\qqq
 \left(\begin{array}{c}0\\1\end{array}\right)=\frac1{\t_+-\t_-}(e_+-e_-).
$$
Then we get $\t_\pm=e^{\pm i\x}$ and $\cos\x={r\/2}$. Assume that
$\t_+\ne\t_-$, then
$$
 \left(\begin{array}{c}\vp_{2k}\\\vp_{2k+1}\end{array}\right)=
 \frac1{\t_+-\t_-}M_1^{k}(e_+-e_-)=\frac1{\t_+-\t_-}(\t_+^{k}e_+-\t_-^{k}e_-),
$$
which yields
\begin{multline}
\lb{fs1}
\vp_{2k}={\l\/a}\frac{\t_+^{k}-\t_-^{k}}{\t_+-\t_-}={\l\/a}{\sin
k\x\/\sin\vk}, \qqq  \ \cos\x={r\/2}\\
  \vp_{2k+1}={1\/a}{\t_+^{k}-\t_-^{k}\/\t_+-\t_-}+
{\t_+^{k+1}-\t_-^{k+1}\/\t_+-\t_-}={\sin
k\x\/a\sin\vk}+{\sin(k+1)\x\/\sin\x}.
\end{multline}
All our functions are analytic, then  \er{fs1} holds true for all
 $\l$, except finite numbers of $\l$. Moreover,  $\vp_n$ is a polynomial,
 then we obtain \er{fs1} for all $\l\in \C$.

In order to determine $\l_k^0(a)$ we calculate all zeros of the
function $\vp_{p+1}(\cdot, a)=0$ and \er{fs1} yields \er{T2-1},
i.e., $\l_{0}^0(a)=0$ and $\l_{-k}^0(a)=-\l_k^0(a), k\in\N_N$, where
\[
\lb{l0k}
 \l_k^0(a)=|a^2-2ac_k+1|^{1\/2}
 =|(a-c_k)^2+s_k^2|^{1\/2}\qq  all\ a\in [0,2].
\]
and $c_k=\cos{k\pi\/N+1}$, $s_k=\sin{k\pi\/N+1}$. The proof of
\er{T2-2}-\er{T2-3} follows from the analysis of the function
$\l_k^0(a)$, see \er{l0k}. This function is monotonic on the
intervals $[0,c_k], [c_k,2]$ and has the extremum at $a=c_k$ such
that $\l_k^0(c_k)=s_k$. Note that $c_k\in[0,1]$ only for $|k|\le
{N+1\/2}$.
 \BBox

\section {Proof of main Theorems }
\setcounter{equation}{0}

 In this section we consider the spectrum of $H$ as $v\to 0$.
 Let $\pa_a={\pa\/\pa a}$ and $\pa_j={\pa\/\pa v_j}$.

{\no\bf Proof of Theorem \ref{T3}.} i) Consider the case $k>0$, the
proof for $k<0$ is similar. It is well known (see e.g. \cite{K2})
that each eigenvalue $\l_k(a,v)$ is a real analytic function of
$(a,v)\in (0,2]\ts \R^p$ and in particular $\l_0(\cdot,\cdot)$ is
analytic in the neighborhood of $(0,0)$.  Identity \er{T2-1} gives
that $\pa_a^2\l_k(a,0)>0$, $a\in[\d,2]$ for some small $\d>0$. Then
the analyticity of $\l_k(a,v)$ in $(a,v)$ implies
$\pa_a^2\l_k(a,v)>0$, for all $a\in[\d,2]$ and for sufficiently
small $v$.

ii) {\bf Consider the case $|k|<{N+1\/2}$}. We have
\[
 \lb{t0000}
 \min_{a\in[0,2]}\l_k(a,0)=\l_k(c_k,0)=s_k,\ \ \ s_k<\l_k(a,0),\
 a\in[0,2]\sm c_k.
\]
Using $\pa_a\l_k(c_k,0)=0$ and $\pa_a^2\l_k(c_k,0)>0$ and the
Implicit Function Theorem we deduce that there exists an unique
function $z_k(v)$ (analytic in $\|v\|\le \d$ for some $\d>0$) such
that $\pa_a\l_k(z_k(v),v)=0$ and $\pa_a^2\l_k(z_k(v),v)>0$ for all
$\|v\|\le \d$ and $z_k(0)=c_k$. Finally we have
$$
 \l_k^-(v)=\min_{a\in[0,2]}\l_k(a,v)=\l_k(z_k(v),v),
$$
\[
 \lb{t0002}
 \pa_j\l_k^-(0)=\pa_a\l_k(a_k(0),0)\pa_ja_k(0)+\pa_j\l_k(a_k(0),0)
 =\pa_j\l_k(c_k,0).
\]
Using i) and the perturbation theory,  we obtain
\[
 \lb{t0004}
 \pa_j\l_k(c_k,0)=
 \frac{\vp_j^2(s_k,c_k,0)}{\sum_{n=1}^{p}\vp_n^2(s_k,c_k,0)}.
\]
Note that \er{t0001} and \er{fs1} gives
$$
 2\cos\x_k=\frac{s_k^2}{c_k}-c_k-\frac1{c_k}=-2c_k,\ \
 \x_k=\pi-\frac{k\pi}{N+1},
$$
\begin{multline}
\lb{t0003}
 \vp_{2n}(s_k,c_k,0)={s_k\sin{n\vk_k}\/c_k\sin\vk_k}={(-1)^{n+1}s_{nk}\/c_k},\\
 \ \ \vp_{2n-1}(s_k,c_k,0)={\sin n\vk_k\/c_k\sin\vk_k}
 +{\sin n\vk_k\/\sin\vk_k}=(-1)^{n}{s_{(n-1)k}\/c_ks_k}+
 (-1)^{n-1}{s_{nk}\/s_k}\\= (-1)^{n}{s_{nk}c_k-c_{nk}s_k\/c_ks_k}
 +(-1)^{n-1}{s_{nk}\/s_k}=(-1)^{n-1}{c_{nk}\/c_k}.
\end{multline}
Thus identities \er{t0002}-\er{t0003} imply
$$
\sum_{n=1}^{p}\vp_n^2(s_k,c_k,0)=\frac{N+1}{c^2_k},\qqq
 \pa_j\l_k^-(0)=\ca\frac{s^2_{nk}}{N+1}, & j=2n \\
 \frac{c^2_{nk}}{N+1}, & j=2n-1 \ac,
$$
and the Taylor formula yields  \er{T3-1}.

{\bf Consider $\l^-_k(v)$ for the case $|k|\ge {N+1\/2}$}.
Asymptotics \er{T3-2} follow from \er{T2-1}, since  $\l_k(a,0)$ is
monotonic in $a\in[0,2]$ and $\l_k(0,0)=\sign k$.

{\bf Consider the asymptotics of $\l^+_k(v)$} for the case $k\ge 1$,
the proof for $k<0$ is similar. We have that
\[
 \lb{t0010}
 \l_k^+(v)=\max_{a\in[0,2]}\l_k(a,v)=\l_k(2,v)
\]
for sufficiently small $v$, since at $v=0$ it follows from \er{T2-1}
and for small $v$ it follows from analyticity of  $\l_k(a,v)$ on
$(0,2]\ts \R^p$. Then the perturbation theory yields
\[
 \lb{t0011}
 \pa_j\l_k(2,0)=
 \frac{\vp_j^2(\m,2,0)}{\sum_{n=1}^{p}\vp_n^2(\m,2,0)},\
 \ \m= \l_k^{0+}=\l_k(2,0)=(5-4c_k)^{\frac12}.
\]
We will calculate $\vp_n(\m,2,0)$. Using again  \er{t0001} and
\er{fs1} we obtain
$$
 2\cos\x_k={\m^2\/2}-2-{1\/2}=-2c_k,\ \
 \x_k=\pi-\frac{k\pi}{N+1},
$$
\[ \lb{t0012}
 \vp_{2n}^2(\m,2,0)={\m^2\sin^2n\x_k\/4\sin^2\x_k}=
 {\m^2s_{kn}^2\/4s_k^2}={\m^2\/8s_k^2}(1-c_{2nk}),\ \
 \
 \sum_{n=1}^{N}\vp_{2n}^2(\m,2,0)=\frac{\m^2(N+1)}{8s_k^2},
\]
and
\begin{multline}
 \lb{t0013}
 \vp_{2n+1}^2(\m,2,0)=\lt({\sin
n\x_k\/2\sin\x_k}+
{\sin(n+1)\x_k\/\sin\x_k}\rt)^2=\lt({s_{nk}\/2s_k}-{s_{(n+1)k}\/s_k}\rt)^2
\\
 ={s_{nk}^2\/4s_k^2}+{s_{(n+1)k}^2\/s_k^2}-{s_{nk}s_{(n+1)k}\/s_k^2}=
{1-c_{2nk}+4(1-c_{2(n+1)k})-4(c_k-c_{(2n+1)k})\/8s_k^2}
\\
 ={\m^2\/8s_k^2}-{c_{2nk}+4c_{2nk+2k}-4c_{2nk+k}\/8s_k^2},\
 \qqq
 \sum_{n=0}^N\vp_{2n+1}^2(\m,2,0)={\m^2(N+1)\/8s_k^2},
\end{multline}
\[
\lb{t0014}
 \sum_{n=1}^{2N+1}\vp_{n}^2(\m,2,0)=
 \sum_{n=0}^N\vp_{2n+1}^2(\m,2,0)+\sum_{n=1}^{N}\vp_{2n}^2(\m,2,0)=
 \frac{\m^2(N+1)}{4s_k^2}
\]
where the identity $\sum_{n=0}^Nc_{nk+\a}=\sum_{n=0}^Ns_{nk+\a}=0$
was used for all $1\le |k|\le N$ and $\a\in\R$. Then substituting
\er{t0012}-\er{t0014} into \er{t0011} and using \er{t0010}, the
Taylor formula gives \er{T3-3}. \BBox

 We consider now the "small spectral band" $\s_0(v)$ of $H$ as $v\to 0$.

\no {\bf Proof of Theorem \ref{T4a}.} i) {\bf Sufficiency.} If
$v_1=v_3=..=v_{2N+1}=0$, then \er{mm1} gives
$$
 T_k=T_k(0,a,v)={1\/a}\ma -1& 0\\
 v_{2k}& -a^2
 \am, \ \ and\ \ T_{N+1}T_N..T_1=\ma .. & 0\\ ..&..\am,
$$
which yields $\vp_{p+1}(0,a,v)=0$ for all $a\in(0,2]$. Then $\l=0$
is an eigenvalue of each $J_a, a\in [0,2]$.

{\no\bf Necessity.} Without loss of generality we assume that
$v_1=0$. Suppose $\l=0$ is an eigenvalue of $H$. Then
$\vp_{p+1}(0,a,v)=0$ for infinity numbers of $a\in(0,2]$. But
$a^{N+1}\vp_{p+1}$ is a polynomial of $a$, then we obtain
$\vp_{p+1}(0,a,v)\ev0$ for any $a$. We have
$$
 T_k=T_k(0,a,v)={1\/a}\ma -1& -v_{2k-1}\\
 v_{2k}& v_{2k}v_{2k-1}-a^2
 \am.
$$
Thus due to \er{mm1} simple calculations gives
$$
 M_k(0,a,v)=\frac1{a^k}\ma .. & O(a^{2k-1}) \\ .. &
 (-1)^ka^{2k}+O(a^{2k-1}).
 \am.
$$
Then
$$
 M_{N+1}(0,a,v)=\frac1{a^{N+1}}\ma -1& -v_{p}\\
 v_{p+1}& v_{p+1}v_{p}-a^2
 \am\ma .. & O(a^{2N-1}) \\ .. & (-1)^Na^{2N}+O(a^{2N-1})
 \am
$$
$$
 =\frac1{a^{N+1}}\ma .. & (-1)^{k+1}v_{p}a^{2N}+O(a^{2N-1}) \\ .. &
 ..\am,
$$
which yields
$a^{N+1}\vp_{p+1}(0,a,v)=(-1)^{N+1}v_{p}a^{2N}+O(a^{2N-1})$. Using
identity $\vp_{p+1}(0,a,v)=0$ for all $a$, we obtain $v_{p}=0$ and
then
$$
 M_{N+1}=\frac1{a}\ma -1& 0\\
 v_{2N+2}& -a^2 \am M_N=\frac1{a^{N+1}}\ma .. &
 -a^{-1}\vp_{2N}(0,a,v)\\..&..\am,
$$
which yields $\vp_{2N}(0,a,v)=-a\vp_{2N+2}(0,a,v)=0$ for all $a$.
Repeat arguments as above for the function $\vp_{2N}(0,a,v)=0$,
$a\in\R$ we obtain $v_{2N-1}=0$ and so on.

ii) Now we describe eigenfunction
$f=(f_{n,k})_{\Z\ts\N_{p}}\in\ell^2(\G)$ such that $Hf=0$. Recall
that the statement i) yields that each $v_{2k+1}=0, k\in\N^0_{N}$.
Define $\p_k=(f_{n,k})_{n\in\Z}$, $k\in\N_{p}$. Then, using \er{Km},
we obtain that $\p=(\p_k)_1^{2N+1}$ satisfies $K\p=0$ and \er{Km}
gives the first equation
$$
 v_1\p_1+(I+S^{-1})\p_2=(I+S^{-1})\p_2=0.
$$
Any $\p_1\in\ell^2(\Z)$ and $\p_2=0$ are all solutions of this
equation. Substituting $\p_2=0$ into the second equation in \er{Km},
we obtain
$$
 (I+S)\p_1+v_2\p_2+\p_3=(I+S)\p_1+\p_3=\l_0\p_2=0,
$$
which yields $\p_3=(-I-S)\p_1$ and so on.

iii)  Let $\wt f=\sum_{m}\wh f_m\p^m$. The definition of $\wt f$
yields $h=(\wt f_{n,k})_{n\in \Z}=(f_{n,k})_{n\in \Z}$ and the
statement ii) gives $\wt f=f$. The similar arguments imply that the
mapping $f\to (\wh f_{m})_{m\in\Z}$ is a linear isomorphism between
$\cH$ and $\ell^2(\Z)$. \BBox

\no {\bf Proof of Theorem \ref{T4}.} i) Consider the eigenvalue
$\l_0(a,v)$, which is a zero of $R(\l,a,v)=a^{N+1}\vp_{p+1}(\l,a,v)$
and $\l_0(a,v)$ is a perturbation of $\l_0^0(a)=\l_0(a,0)$. Note
that $R$ is a polynomial and $R(\l,a,v)=\l^{p}+O(\l^{2N})$ as $\l\to
\iy$.

Consider $R(\l,a,0)$. Recall that $\l_0^0(a)=0$ is a simple zero of
$R(\l,a,0)$ for all $a\in[0,2]$. The function $ R(\l,a,v)$ is a
polynomial of degree $p$ in $\l$, whose coefficients are polynomials
in the components of $a,v$. Therefore, its roots $\l_n, n\in \Z_N$
are continuous functions of $(a,v)\in [0,2]\ts \R^p$. Moreover,
these roots are simple, then they are real analytic on $[0,2]\ts
\R^p$.

Using Theorem \ref{T3} i), we obtain
$$
 \l_{0}(a,v)=\l_{0}(a,0)+\sum_{k=1}^{p}\pa_k\l_{0}(a,0)v_k+O(\|v\|^2),
$$
uniformly in $|a|+|v|\le \d$ for some small $\d>0$. Using \er{fs1}
and \er{t0001} we deduce that
$$
 \vp_{2k}(0,a,0)=0,\ \ \vp_{2k+1}(0,a,0)=(-a)^k,\qq
 \sum_{j=1}^{p}\vp_k^2(0,a,0)=\sum_{k=0}^Na^{2k}=A,
$$
which yields
$$
 \pa_{2k}\l_{0}(a,0)={\vp_{2k}^2(0,a,0)\/A}=0,\qqq
 \pa_{2k+1}\l_{0}(a,0)={\vp_{2k+1}^2(0,a,0)\/A}=
 \frac{a^{2k}(a^2-1)}{a^{2N+2}-1}.
$$
Thus
$$
 \l_{0}(a,v)=F(a,v)+O(\|v\|^2),\qqq
 F(a,v)={\sum_{k=0}^Nv_{2k+1}a^{2k}\/\sum_{k=0}^Na^{2k}}.
$$
which yields \er{T5-1}, \er{T5-1a}.

ii) Let $v\to0$. Using   $F(0,v)=0$ and $
F(2,v)={3\/4^{N+1}-1}\sum_{k=1}^N4^{k}v_{2k+1} $ and \er{T5-1},
\er{T5-1a} and Lemma \ref{Tp}, we obtain \er{T5-2}. Moreover, If
$v_1<v_3$ and $v_k\le C v_1$ then $\l_0(a,v)$  is strictly
increasing for $a\in[0,2]$, see Lemma \ref{L7a}. \BBox

\begin{lemma}
\lb{Tp}
Let $v_{1}\le v_{3}\le..\le v_{2N+1}$ and  $v_{1}<
v_{2N+1}$. Then the function $F(\cdot,v)$ is strictly increasing on
$[0,2]$ and given by
\[ \lb{Tp-1}
 F(a,v)=v_{2N+1}-\sum_{k=1}^N(v_{2k+1}-v_{2k-1})\frac{a^{2k}-1}{a^{2(N+1)}-1}.
\]
\end{lemma}
{\no\bf Proof.} We have that
$$
 F(a,v)=
 \frac{(\sum_{0}^Nv_{2k+1}a^{2k})(a^2-1)}{a^{2(N+1)}-1}
 =\frac{v_{2N+1}a^{2(N+1)}+\sum_{1}^N(v_{2k-1}-v_{2k+1})a^{2k}-v_1}{a^{2(N+1)}-1}
$$
$$
 =\frac{v_{2N+1}(a^{2(N+1)}-1)+\sum_{1}^N(v_{2k-1}-v_{2k+1})(a^{2k}-1)}{a^{2(N+1)}-1}
 =v_{2N+1}-\sum_{k=1}^N(v_{2k+1}-v_{2k-1})\frac{a^{2k}-1}{a^{2(N+1)}-1}.
$$
Then the function $F(\cdot,v)$ is strongly increasing on $[0,2]$,
since the function $-{x-1\/x^r-1}$  is strongly increasing in $x\ge
0$ for any $r>1$ and then $-\frac{a^{2k}-1}{a^{2(N+1)}-1}$ is
strictly increasing for $a>0$.
\BBox

\begin{lemma}
\lb{L7a}
The function $\m(z,v)=\l_0(\sqrt z,v)$ is real analytic in
$$
 (z,v)\in\O_{\d}=\O^1_{\d}\ts \cB_{\d},\ \ \O^1_\d=\{z\in\C:\ \dist(z,[0,4])<\d\},\ \
 \cB_\d=\{v\in\R^p:\ \|v\|<\d\},
$$
for some $\d>0$.  Moreover, the function  $\wt
F(z,v)=\m(z,v)-F(\sqrt z,v)$ satisfies
\[
 \lb{L32-2}
 |\wt F(z,v)|+|{\pa_z}\wt F(z,v)|\le C_0\|v\|^2,\ \ (z,v)\in\O_{\d},
\]
for some constant $C_0$. Furthermore, if $0=v_1<v_3\le..\le
v_p<cv_3$ for some $c>0$, then $\pa_z\m(z,v)>0$ for all
$(z,v)\in[0,4]\ts \cB_\d$  for sufficiently small $\d>0$.

\end{lemma}
\no {\bf Proof.} The identities \er{mm1},\er{mo} imply
$\vp_{p+1}(\l,-a,v)=(-1)^{N+1}\vp_{p+1}(\l,a,v)$. Then  we obtain
$\l_0(a,v)=\l_0(-a,v)$ for sufficiently small $v$, since
$\vp_{p+1}(\l_0(a,v),a,v)=0$ and $\l_0(a,v)$ is a simple root of the
polynomial $\vp_{p+1}(\l,a,v)$ for sufficiently small $v$. Then, the
Implicit Function Theorem gives that $\l_0(a,v)$ is real analytic in
$\O_{3\d}$ and the new function $\m(z,v)=\l_0(\sqrt z,v), z=a^2$ is
real analytic in $\O_{3\d}$, since $\l_0(a,v)=\l_0(-a,v)$.

This result  and \er{T5-1} yield $|\wt F(z,v)|\le C_1\|v\|^2$ for
all $(z,v)\in\O_{3\d}$  and some $C_1$.

If $(z_0,v)\in\O_{\d}$, then $\m(z,v)$ is analytic in $z\in
\{|z-z_0|<\d\}$ and the Cauchy integral implies
$$
\pa_z \wt F(z_0,v)={1\/2\pi i}\int_{|z-z_0|=\d}\frac{\wt
F(z,v)}{(z-z_0)^2}dz,\qqq
 |\pa_z \wt F(z_0,v)|\le {C_1\/\d}\|v\|^2,
$$
since $|\wt F(z,v)|\le C_1\|v\|^2$ for all $(z,v)\in\O_{3\d}$.

Let $f_k(z)={1-z^{k}\/z^{N+1}-1},\ \ z>0,$ for fix $k$. This
function is strongly increasing for $z\ge0$ and
$A=\min_{z\in[0,4]}(f_1)'(z)>0$. Then using \er{Tp-1}, \er{L32-2},
we obtain
$$
 \pa_z\m(z,v)=\pa_z G(z,v)+\pa_z\wt
 F(z,v)=\sum_{k=1}^N(v_{2k+1}-v_{2k-1})(f_k)'(z)+\pa_z\wt
 F(z,v)\ge v_3A+O(v_3^2)
$$
which yields $\pa_z\m(z,v)>0$   for all sufficiently small $v_3$.
\BBox

\no {\bf Proof of Theorem \ref{T5}.} Using \er{Jk1} we rewrite
$J_a(tv)$ in the form
$$
 J_a(tv)=t(B+t^{-1}J_a^0)=t(B+\ve J_a^0),\
                                  \ B=\diag(v_j)_{1}^{p},\ as \
                                   \ \ve={1\/t}\to 0.
$$
Then the perturbation theory (see  Sect. XII, 1, \cite{RS})  for
$B+\ve J_a^0$ gives
\begin{multline}
\lb{T5.b}
 {\l_{k-N-1}(a,tv)}=t(v_k-\a_k\ve^2+O(\ve^3)),\ \
 \a_k=\sum_{j\in \N_p\sm\{k\}}{u_{k,j}^2\/v_j-v_k},\ \
 u_{j,k}=(e_j^0,J_a^0e_k^0),
\end{multline}
for any $k\in\N_p$, here  $Be_j^0=v_je_j^0$ and the vector
$e_j^0=(\d_{j,n})_{n=1}^{p}\in \C^{p}$. The matrix
$J_a^0=\{u_{k,j}\}$ is given by \er{ja0}, where
\[ \lb{T5.b1}
 u_{j,k}\ev u_{j,k}(a)=\ca0, & |k-j|\ne1,\\
 a, & j=k+1,\ j\in2\N,\ or\ j=k-1,\ k\in2\N
 \\ 1, & other\ cases \ac\ .
\]
Then, each function $u_{j,k}(a)$ is strongly increasing (or
constant)  in $a\in[0,2]$. Using similar arguments as in Theorem
\ref{T4} and Lemma \ref{L7a}, we deduce that each function
$\l_k(a,tv)$, $k\in\Z_N\sm\{N\}$ is strongly increasing  in
$a\in[0,2]$ and for sufficiently large $|t|$. Note that if $k=N$,
then $\l_N(a,tv)=tv_p+{1\/t(v_p-v_{p-1})}+O(t^{-2})$ and our simple
arguments do not give that $\l_N(a,tv)$ is monotonic in $a\in[0,2]$.
Using these arguments and $\s_k(t)\cap\s_n(t)=\es$, $k\ne n$ (see
\er{T5.b}), we deduce that $\s_k(t), k\ne N$ has multiplicity $2$
for sufficiently large $t$. The endpoints of $\s_k(t)$, $ k\ne N$
are $\l_k(0,tv)$ and $\l_k(2,tv)$, then \er{T5.b}, \er{T5.b1}  yield
\er{T5.1}, \er{T5.3}.
\BBox

 \no {\bf Acknowledgments.}
 \small
E. K. was partly supported by DFG project BR691/23-1. The various
parts of this paper were written at the Erwin Schr\"odinger
Institute for Mathematical Physics, Vienna, E.K. is grateful to the
Institutes for the hospitality.  The some parts of this paper was
written at the Math. Institute of Humboldt Univ., Berlin; A. K. is
grateful to the Institute for the hospitality.

\end{document}